\documentclass[12pt]{article}
\title
{The critical bias for the Hamiltonicity game is $n/\ln n$}
\author{Michael Krivelevich
\thanks{School of Mathematical Sciences, Raymond and Beverly
Sackler Faculty of Exact Sciences, Tel Aviv University, Tel Aviv,
69978, Israel. Email: krivelev@post.tau.ac.il. Research supported in
part by a USA-Israel BSF grant, by a grant from the Israel Science
Foundation and by a Pazy Memorial Award.}
 }

\usepackage{amsmath, amsfonts}

\begin{document}
\bibliographystyle{plain}
\maketitle
\newtheorem{theorem}{Theorem}
\newtheorem{claim}{Claim}
\newtheorem{lemma}{Lemma}
\newtheorem{propos}{Proposition}
\newtheorem{conjecture}{Conjecture}
\newtheorem{problem}{Problem}
\newtheorem{corol}{Corollary}
\newcommand{\Proof}{\noindent{\bf Proof.}\ \ }
\newcommand{\Remarks}{\noindent{\bf Remarks:}\ \ }
\newcommand{\whp}{{\bf whp}}
\newcommand{\prob}{probability}
\newcommand{\rn}{random}
\newcommand{\rv}{random variable}
\newcommand{\hpg}{hypergraph}
\newcommand{\hpgs}{hypergraphs}
\newcommand{\subhpg}{subhypergraph}
\newcommand{\subhpgs}{subhypergraphs}
\newcommand{\bH}{{\bf H}}
\newcommand{\cH}{{\cal H}}
\newcommand{\cT}{{\cal T}}
\newcommand{\cF}{{\cal F}}
\newcommand{\cD}{{\cal D}}
\newcommand{\cC}{{\cal C}}

\begin{abstract}
We prove that in the biased $(1:b)$ Hamiltonicity Maker-Breaker
game, played on the edges of the complete graph $K_n$, Maker has a
winning strategy for $b(n)\le
\left(1-\frac{30}{\ln^{1/4}n}\right)\frac{n}{\ln n}$, for all large
enough $n$.
\end{abstract}

\section{Introduction}
A {\em Maker-Breaker game} is a triple $(H,a,b)$, where $H=(V,E)$ is
a hypergraph with vertex set $V$, called the board of the game, and
edge set $E$, a family of subsets of $V$ called winning sets. The
parameters $a$ and $b$ are positive integers, related to the so
called game bias. The game is played between two players, called
Maker and Breaker, who change turns occupying previously unclaimed
elements of $V$; Maker claims $a$ elements in his turn, Breaker
answers by claiming $b$ elements. We assume that Breaker moves
first. The game ends when all board elements have been claimed by
either of the players. (In the very last move, if the board does not
contain enough elements to claim for the player whose turn is now,
that player claims all remaining elements of the board.) Maker wins
if and only if he has occupied one of the winning sets $e\in E$ by
the end of the game, Breaker wins otherwise, i.e., if he manages to
occupy at least one element of (``to break into") every winning set
by the end of the game. The most basic case is when $a=b=1$, which
is the so called {\em unbiased} game. Here we will be concerned with
$1:b$ games.

It is quite easy to see that Maker-Breaker games are {\em bias
monotone}. This is to say that if the game $(H,1,b)$ is Maker's win,
then $(H,1,b')$ is Maker's win as well for every integer $b'<b$.
This allows to define the {\em critical bias} of the game $H$, which
is the maximum possible value of bias $b$ for which Maker still wins
the $1:b$ game played on $H$ (if the 1:1 game is Breaker's win, we
say that the critical bias in this case is zero).

We refer the reader to a recent monograph \cite{Beck-book} of Beck
for extensive background on positional games in general and on
Maker-Breaker games in particular.


The subject of this paper is the {\em Hamiltonicity game} played on
the edge set of the complete graph $K_n$. In this game, players take
turns in claiming unoccupied edges of $K_n$. Maker's aim is to
construct a Hamilton cycle, and thus the family of winning sets
coincides with the family of (the edge sets of) graphs on $n$
vertices containing a Hamilton cycle. The research on biased
Hamiltonicity games has a long and illustrious history. Already in
the very first paper about biased Maker-Breaker games back in 1978,
Chv\'atal and Erd\H{o}s \cite{CE} treated the unbiased Hamiltonicity
game and showed that Maker wins it for every sufficiently large $n$.
(Chv\'atal and Erd\H{o}s showed in fact that Maker wins within $2n$
rounds. Later the minimum number of steps required for Maker to win
this game was shown to be at most $n+2$ by Hefetz et al. in
\cite{HKSS}, and finally the optimal $n+1$ by Hefetz and Stich
\cite{HS}.) We would like to mention that Chv\'atal and Erd\H{o}s
also proved in their paper that for $b(n)\ge (1+\epsilon)n/\ln n$,
where $\epsilon>0$ is an arbitrary small constant, Breaker can
isolate a vertex in the $1:b$ game played on $K_n$, i.e., to claim
all $n-1$ edges incident to it.

Chv\'atal and Erd\H os conjectured that there is a function $b(n)$
tending to infinity such that Maker can still build a Hamilton cycle
if he plays against bias $b(n)$. Their conjecture was verified by
Bollob\'as and Papaioannou \cite{BP} who proved that Maker is able
to build a Hamilton cycle even if Breaker's bias is as large as
$\frac{c\ln n}{\ln \ln n}$, for some constant $c>0$. Beck improved
greatly on this \cite{Beck85} and showed that Maker wins the
Hamiltonicity game provided Breaker's bias is at most
$\left(\frac{\ln 2}{27}-o(1)\right)\frac{n}{\ln n}$. In view of the
above mentioned Chv\'atal-Erd\H{o}s theorem about isolating a
vertex,  Beck's result established that the order of magnitude of
the critical bias in the Hamiltonicity game is $n/\log n$.
Krivelevich and Szab\'o \cite{KS} improved upon Beck's result and
showed that the critical bias $b(n)$ for the Hamiltonicity game is
at least $(\ln 2-o(1))n/\ln n$.

In a relevant development, Gebauer and Szab\'o showed recently in
\cite{GS} that the critical bias for the connectivity game on $K_n$
(where Maker wins if and only if he creates a spanning tree from his
edges by the end of the game) is asymptotically equal to $n/\ln n$.
We will rely extensively on some of their results and approaches
here.

It was widely believed that the critical bias for the Hamiltonicity
game on $K_n$ is asymptotically equal to $n/\ln n$ as well. This
conjecture has even attained the (somewhat dubious) honor to be
stated as one of the most ``humiliating open problems" of the
subject by Beck in his book \cite{Beck-book} (see Chapter 49 there).

\section{The result}
In this paper we resolve the above stated conjecture. Here is our
result:
\begin{theorem}\label{th1}
Maker has a strategy to win the $(1:b)$ Hamiltonicity game played on
the edge set of the complete graph $K_n$ on $n$ vertices in at most
$14n$ moves, for every
$b\le\left(1-\frac{30}{\ln^{1/4}n}\right)\frac{n}{\ln n}$, for all
large enough $n$.
\end{theorem}

The constants and the error term expression of the above theorem are
clearly not optimal and can be improved somewhat by a more careful
implementation of our arguments. We find however little reason to
pursue this goal.

\section{Notation}
Our basic notation is quite standard and follows closely that of
most of graph theory books. In particular, for a graph $G=(V,E)$ and
a vertex subset $U\subset V$, we denote by $N_G(U)$ the external
neighborhood of $U$ in $G$, i.e., $N_G(U)=\{v\in V\setminus U:
v\mbox{ has a neighbor in } U\}$. We systematically omit rounding
signs for the sake of clarity of presentation. The underlying
parameter $n$ is assumed to be large enough where necessary.

Let
\begin{eqnarray*}
\delta_0 = \delta_0(n) &=& \frac{6}{\ln^{1/2}n}\,,\\
\delta  = \delta(n)  &=& \frac{15}{\ln^{1/4}n}\,,\\
\epsilon = \epsilon(n) &=& \frac{30}{\ln^{1/4}n}\,,\\
k_0 = k_0(n) &=& \delta_0n = \frac{6n}{\ln^{1/2}n}\,.
\end{eqnarray*}

For a positive integer $k$, a graph $G=(V,E)$ is a {\em
$k$-expander} if $|N_G(U)|\ge 2 |U|$ for every subset $U\subset V$
of at most $k$ vertices.

Given a graph $G$, a non-edge $e=(u,v)$ of $G$ is called a {\em
booster} if adding $e$ to $G$ creates a graph $G'$, which is
Hamiltonian or whose maximum path is longer than that of $G$.
Boosters advance a graph towards Hamiltonicity when added; adding
sequentially $n$ boosters clearly brings a graph to Hamiltonicity.

\section{Tools}

The following lemma, that can be traced back to a seminal work of
P\'osa \cite{Pos}, is used quite frequently in papers on
Hamiltonicity and extremal problems involving paths and cycles.

\begin{lemma}\label{le1}
Let $G$ be a connected non-Hamiltonian $k$-expander. Then at least
$(k+1)^2/2$ non-edges of $G$ are boosters.
\end{lemma}

\Proof See, e.g., Lemma 8.5 of \cite{Bol-book} or Corollary 2.10 of
\cite{KLS}. \hfill $\Box$

\medskip

Although $k$-expanders are not necessarily connected, their
connected components are guaranteed to be of a relatively large
size, as shown in the following easy lemma.

\begin{lemma}\label{le2}
Let $G=(V,E)$ be a $k$-expander. Then every connected component of
$G$ has size at least $3k$.
\end{lemma}

\Proof If not, let $V_0$ be the vertex set of a connected component
of $G$ of size less than $3k$. Choose an arbitrary subset
$U\subseteq V_0$ of cardinality $|U|=\min\{|V_0|,k\}$, clearly
$|U|>|V_0|/3$. Since $G$ is a $k$-expander, it follows that
$|N_G(U)|\ge 2|U|$. On the other hand, $N_G(U)\subseteq V_0$,
implying $|V_0|\ge |U|+|N_G(U)|\ge 3|U|$ -- a contradiction. \hfill
$\Box$

\medskip

Now we can describe the main tool of our proof, a recent result of
Gebauer and Szab\'o, who analyzed in \cite{GS} the biased minimum
degree game. For our goals, it will suffice to specialize their
analysis to the game where Maker's goal is to reach a graph of
minimum degree at least 12. Maker's strategy employed by Gebauer and
Szab\'o is very natural:

\medskip
\begin{center}
\fbox{\parbox{6in} {{\bf Strategy $S$:}\\
{
 As long as there is a vertex of degree less than 12 in Maker's
graph, Maker chooses a vertex $v$ of minimum degree in his graph
(breaking ties arbitrarily) and claims an {\bf arbitrary} unclaimed
edge $e$ containing $v$. } }}
\end{center}
\medskip

If Maker claims an edge $e$ due
to a vertex $v$ in the above strategy, we say that $e$ {\em is
chosen by} $v$. Gebauer and Szab\'o proved the following statement
about it.

\begin{theorem} (\cite{GS}, Theorem 1.2):
In a $\big(1:\frac{(1-\epsilon)n}{\ln n}\big)$-game played on the
edge set of the complete graph $K_n$ on $n$ vertices, strategy $S$
guarantees Maker minimum degree at least 12 in his graph.
\end{theorem}

In our argument we will need more than the above statement -- it
will be essential for us that Maker is able, for every vertex $v$ of
the graph, to reach degree at least 12 at $v$ when a substantial
part of the edges at $v$ is still unclaimed. Fortunately, the proof
of Gebauer and Szab\'o gives this as well, as stated in the lemma
below.

\begin{lemma}\label{le-GS}
In a $\big(1:\frac{(1-\epsilon)n}{\ln n}\big)$-game played on the
edge set of the complete graph $K_n$ on $n$ vertices, strategy $S$
guarantees that for every vertex $v\in[n]$ Maker has at least 12
edges incident to $v$ before Breaker accumulates at least
$(1-\delta)n$ edges at $v$.
\end{lemma}

The proof of Lemma \ref{le-GS} is a straightforward modification of
the proof of Theorem 1.2 of \cite{GS}. More specifically, in their
argument one just needs to notice that in the current setting the
danger of the last vertex $v_g$ before Breaker's last move is now at
least $(1-\delta)n-12-b$, and then to check that in the relevant
calculations the danger of the original set $I_{g-1}$ before the
game started still comes out positive. We refer the reader to
\cite{GS} for further details.

\section{The proof}
In this section we prove our main result, Theorem \ref{th1}. Maker's
strategy is composed of three stages. At the first stage, he creates
a $k_0$-expander in a linear number of moves. At the second stage,
Maker makes sure his graph is connected in at most $O(n/k_0)$ moves.
Finally, he turns his graph into a Hamiltonian one, using at most
$n$ further moves.

\medskip
{\bf Stage 1 -- creating an expander.}

Let us go back to the Gebauer-Szab\'o winning strategy $S$ for the
minimum degree 12 game. As it turns out, this strategy not only
guarantees minimum degree 12 or more in Maker's graph, but has
enough flexibility in it to allow Maker to pursue an even more
important goal -- that of creating quickly a good expander from its
edges. First observe that as long as the game is played at this
stage, Maker increases by one the degree of a vertex whose current
degree in his graph is still less than 12. Therefore, Maker wins
this game in at most $12n$ moves. More importantly, while describing
strategy $S$, we stressed that at each round Maker is allowed to
choose an edge $e$ incident to its vertex of minimum degree $v$
arbitrarily. We can utilize this freedom of choice by specifying
that Maker claims each time a {\bf random} edge $e$ incident to $v$.
This random choice of Maker allows us to prove that he has a
strategy to create a good expander quickly.
\begin{lemma}\label{le3}
Maker has a strategy to create a $k_0$-expander in at most $12n$
moves.
\end{lemma}

\Proof Maker augments the strategy $S$ described above by choosing
at each round a random edge incident to a vertex of minimum degree
in his graph. Here is his strategy $S'$.

\medskip
\begin{center}
\fbox{\parbox{6in} {{\bf Strategy $S'$:}\\
{
 As long as there is a vertex of degree less than 12 in Maker's
graph, Maker chooses a vertex $v$ of minimum degree in his graph
(breaking ties arbitrarily) and claims a {\bf random} unclaimed edge
$e$ containing $v$. } }}
\end{center}
\medskip

The game lasts till the minimum degree in Maker's graph is at least
12. As we argued above, the game duration does not exceed $12n$.
Since the game analyzed is a perfect information game with no chance
moves, it is enough to prove that Maker's strategy succeeds to
create a $k_0$-expander with positive probability. (We will in fact
prove that his strategy succeeds with probability approaching 1.)

So suppose that Maker's graph is not a $k_0$-expander. Then there is
a subset $A$ of size $|A|=i\le k_0$ in Maker's graph $M$ after the
end of Stage 1 such that $N_M(A)$ is contained in a set $B$ of size
at most $2i-1$. Since the minimum degree in Maker's graph is 12, we
can assume that $i\ge 5$; more importantly, there are at least $6i$
edges of Maker incident to $A$. Consider one such edge $e=(u,v)$ and
assume that $e$ was chosen by $v\in A\cup B$ in the course of the
game. Notice crucially that, by Lemma \ref{le-GS}, when choosing $e$
Breaker's degree at $v$ was at most $(1-\delta)n$, while Maker's
degree at $v$ was at most 11. Therefore at that point of the game,
there were at least $\delta n-12$ unclaimed edges incident to $v$.
The probability that at that point Maker chose an edge at $v$ whose
second endpoint belongs to $A\cup B$ is thus at most $\frac{|A\cup
B|-1}{\delta n-12}$, regardless of the history of the game so far.
It follows that the probability that all these $6i$ edges incident
to $A$ will end up entirely in $A\cup B$ is at most
$\left(\frac{3i-2}{\delta n-12}\right)^{6i}$. Summing over all
relevant values of $i$, we derive that the probability that Maker's
strategy fails to create a $k_0$-expander is at most
\begin{gather*}
\sum_{5\le i\le k_0}{n\choose i}
{{n-i}\choose{2i-1}}\left(\frac{3i-2}{\delta n-12}\right)^{6i} \le
\sum_{5\le i\le k_0} \left[
\frac{en}{i}\,\left(\frac{en}{2i}\right)^2\, \left(\frac{4i}{\delta
n}\right)^6\right]^i\\
=\sum_{5\le i\le k_0}\left[4^5e^3\, \left(\frac{i}{n}\right)^3\,
\frac{1}{\delta^6}\right]^i\ .
\end{gather*}
Denote the $i$-th term of the above sum by $g(i)$. Then for $5\le
i\le \sqrt{n}$ we have $g(i)\le \left(O(1)\ln^{3/2}n\cdot
n^{-3/2}\right)^6= o(1/n)$, while for $\sqrt{n}\le i\le k_0$ we can
estimate $g(i)\le
\left(\frac{4^5e^3\delta_0^3}{\delta^6}\right)^{\sqrt{n}}= o(1/n)$
as well. This implies that Maker's strategy fails with negligible
probability, and thus with positive probability (and in fact almost
surely) he creates a $k_0$-expander in the first $12n$ moves.

\medskip
{\bf Stage 2 -- creating a connected expander.}

If Maker's graph $M$ is not yet connected by the end of Stage 1, he
can turn it easily into such in very few moves. Indeed, $M$ is a
$k_0$-expander and therefore by Lemma \ref{le2} all connected
components of $M$ are of size at least $3k_0$. In the next
$n/(3k_0)-1$ rounds at most, Maker claims an arbitrary edge between
two of its connected components. Observe that there are at least
$9k_0^2=324n^2/\ln n$ edges of the complete graph between any two
such components, and Breaker has at most $(12n+n/(3k_0))\cdot b<
13n^2/\ln n$ edges claimed on the board altogether. Therefore,
Breaker cannot block Maker from achieving his goal. Stage 2 lasts at
most $n/(3k_0)-1<n$ rounds.

\medskip
{\bf Stage 3 -- completing a Hamilton cycle.}

Recall that by the end of Stage 1 Maker has created a
$k_0$-expander. Clearly, his graph at every subsequent round
inherits this expansion property. Also, after Stage 2 Maker's graph
is already connected. But then by Lemma \ref{le1} at any round of
Stage 3 Maker's graph is either already Hamiltonian, or has at least
$k_0^2/2$ boosters. Maker goes on to add a booster after a booster
in the next $n$ rounds at most, till finally he reaches
Hamiltonicity. Breaker is helpless -- he just does not have enough
stuff on the board to block all of Maker's boosters during these
rounds. Indeed, the game lasts altogether at most $12n+n+n=14n$
rounds, during which Breaker puts on the board at most $14n\cdot
b\le 14n^2/\ln n$ edges -- less than $k_0^2/2$ boosters of Maker.
Hence, at any round of Stage 3 there is an available booster with
respect to the current Maker's graph -- which he happily claims.
\hfill $\Box$

\section{Concluding remarks}
We have essentially resolved completely the biased Hamiltonicity
game on the complete graph $K_n$ by proving that the critical bias
$b(n)$ is asymptotic to $n/\ln n$.

The method we employed in our proofs (creating quickly a good
expander first) is quite general and has a clear potential to be
applicable to other biased combinatorial games as well. For example,
it can be used to show that Maker can create a $c$-connected
spanning graph $G$ in the $1:b$ game on $K_n$ for any constant $c$,
or even for a growing function $c=c(n)$, as long as the bias $b(n)$
satisfies $b(n)\le (1-o(1))n/\ln n$. This would provide an
alternative proof of the corresponding results of Gebauer and
Szab\'o \cite{GS} and in fact would strengthen their assertions.

Finally, let us mention that the strategy we used in our argument is
random. It would be very interesting to provide a deterministic
(explicit) Maker's strategy for winning the Hamiltonicity game close
to the critical bias.

\medskip

\noindent{\bf Acknowledgement.} The author wishes to thank Dan
Hefetz for his careful reading of the first draft of the paper and
many helpful remarks.

\end{document}